\documentclass[12pt,a4paper]{article}
\usepackage{amssymb}

\usepackage{graphicx,amssymb,amsfonts,epsfig,amsthm,a4,amsmath,url}
\usepackage[latin1]{inputenc}

\newtheorem{ThmIntro}{Theorem}
\newtheorem{CorIntro}[ThmIntro]{Corollary}

\newtheorem{thm}{Theorem}[section]

\newtheorem{lem}[thm]{Lemma}
\newtheorem{clai}[thm]{Claim}
\newtheorem{prop}[thm]{Proposition}
\theoremstyle{definition}
\newtheorem{defn}[thm]{Definition}

\theoremstyle{remark}
\newtheorem{rem}[thm]{Remark}
\newtheorem{ex}[thm]{Example}

\numberwithin{equation}{section}

\newcommand{\Z}{\mathbf{Z}}

\newcommand{\N}{\mathbf{N}}
\newcommand{\R}{\mathbf{R}}
\newcommand{\C}{\mathbf{C}}
\newcommand{\Q}{\mathbf{Q}}

\newcommand{\D}{\mathbf{D}}

\newcommand{\supp}{\text{Supp}}

\newcommand{\bpr}{\noindent \textbf{Proof}: ~}

\newcommand{\epr}{~$\blacksquare$}

\newcommand{\eps}{\varepsilon}

\title{Vanishing of the first reduced cohomology with values in an $L^p$-representation.}
\author{Romain Tessera}
\date{\today}

\begin{document}

\baselineskip=16pt

\maketitle

\begin{abstract}
We prove that the first reduced cohomology with values in a mixing
$L^p$-representation, $1<p<\infty$, vanishes for a class of
amenable groups including connected amenable Lie groups. In
particular this solves for this class of amenable groups a
conjecture of Gromov saying that every finitely generated amenable
group has no first reduced $\ell^p$-cohomology. As a byproduct, we
prove a conjecture by Pansu. Namely, the first reduced
$L^p$-cohomology on homogeneous, closed at infinity, Riemannian
manifolds vanishes. We also prove that a Gromov hyperbolic
geodesic metric measure space with bounded geometry admitting a
bi-Lipschitz embedded 3-regular tree has non-trivial first reduced
$L^p$-cohomology for large enough $p$. Combining our results with
those of Pansu, we characterize Gromov hyperbolic homogeneous
manifolds: these are the ones having non-zero first reduced
$L^p$-cohomology for some $1<p<\infty.$
\end{abstract}




\section{Introduction}

\subsection{A weak generalization of a result of Delorme.}

In \cite{Del}, Delorme proved the following deep result: every
connected solvable Lie groups has the property that every weakly
mixing\footnote{A unitary representation is called weakly mixing if
it contains no finite dimensional sub-representation.} unitary
representation $\pi$ has trivial first reduced cohomology, i.e.
$\overline{H}^1(G,\pi)\neq 0$. This was recently extended to
connected amenable Lie groups, see \cite[Theorem~3.3]{Mar}, and to a
large class of amenable groups including polycyclic groups by Shalom
\cite{Shal}. Shalom also proves that this property, that he calls
Property $H_{FD}$, is invariant under quasi-isometry between
amenable discrete groups. Property $H_{FD}$ has nice implications in
various contexts. For instance, Shalom shows that an amenable
finitely generated group with Property $H_{FD}$ has a finite index
subgroup with infinite abelianization \cite[Theorem~4.3.1]{Shal}. In
\cite{coteva}, we prove~\cite[Theorem~4.3]{coteva} that an amenable
finitely generated group with Property $H_{FD}$ cannot
quasi-isometrically embed into a Hilbert space unless it is
virtually abelian.

It is interesting and natural to extend the definition of Property
$H_{FD}$ to isometric representations of groups on certain classes
of Banach spaces.

In this paper, we prove that a weak version of Property $H_{FD}$,
also invariant under quasi-isometry, holds for isometric
$L^p$-representations of a large class of amenable groups including
connected amenable Lie groups and polycyclic groups: for
$1<p<\infty,$ every {\it strongly mixing} isometric
$L^p$-representation $\pi$ has trivial first reduced cohomology (see
Section~\ref{mainresultsection} for a precise statement).

\subsection{$L^p$-cohomology.}

The $L^p$-cohomology (for $p$ non necessarily equal to $2$) of a
Riemannian manifold has been introduced by Gol'dshtein, Kuz'minov,
and Shvedov in \cite{GKS}. It has been intensively studied by
Pansu~\cite{Pan,Pansu,Pans} in the context of homogeneous Riemannian
manifolds and by Gromov \cite{Gr} for discrete metric spaces and
groups. The $L^p$-cohomology is invariant under quasi-isometry in
degree one \cite{HS}. But in higher degree, the quasi-isometry
invariance requires some additional properties, like for instance
the uniform contractibility of the space \cite{Gr} (see also
\cite{BP,Pans}). Most authors focus on the first reduced
$L^p$-cohomology since it is easier to compute and already gives a
fine quasi-isometry invariant (used for instance in \cite{B,BP}).
The $\ell^2$-Betti numbers of a finitely generated group,
corresponding to its reduced $\ell^2$-cohomology\footnote{We write
$\ell^p$ when the space is discrete.}, have been extensively studied
in all degrees by authors like Gromov, Cheeger, Gaboriau and many
others. In particular, Cheeger and Gromov proved in \cite{CG} that
the reduced $\ell^2$-cohomology of a finitely generated amenable
group vanishes in all degrees. In \cite{Gr}, Gromov conjectures that
this should also be true for the reduced $\ell^p$-cohomology. For a
large class of finitely generated groups with infinite center, it is
known \cite{Gr,K} that the reduced $\ell^p$-cohomology vanishes in
all degrees, for $1<p<\infty$. The first reduced $\ell^p$-cohomology
for $1<p<\infty$ is known to vanish \cite{BMV,MV} for certain
non-amenable finitely generated groups with ``a lot of
commutativity" (e.g. groups having a non-amenable finitely generated
normal subgroup with infinite centralizer).

\medskip

A consequence of our main result is to prove that the first reduced
$\ell^p$-cohomology, $1<p<\infty$, vanishes for large class of
finitely generated amenable groups, including for instance
polycyclic groups.

On the other hand, it is well known~\cite{Gr} that the first reduced
$\ell^p$-cohomology of a Gromov hyperbolic finitely generated group
is non-zero for $p$ large enough. Although the converse is
false\footnote{In \cite{cotevaBP} for instance, we prove that any
non-amenable discrete subgroup of a semi-simple Lie group of rank
one has non-trivial reduced $L^p$-cohomology for $p$ large enough.
On the other hand, non-cocompact lattices in SO$(3,1)$ are not
hyperbolic. See also \cite{BMV} for other examples.} for finitely
generated groups, we will see that it is true in the context of
connected Lie groups. Namely, a connected Lie group has non-zero
reduced first $L^p$-cohomology for some $1<p<\infty$ if and only if
it is Gromov hyperbolic.

\bigskip

\noindent \textbf{Acknowledgments.} I would like to thank Pierre
Pansu, Marc Bourdon and Herv\'{e} Pajot for valuable discussions
about $L^p$-cohomology. Namely, Marc explained to me how one can
extend a Lipschitz function defined on the boundary
$\partial_{\infty} X$ of a Gromov hyperbolic space $X$ to the
space itself, providing a non-trivial element in $H_p^1(X)$ for
$p$ large enough (see the proof of Theorem~\ref{hyperbolicLpthm}
in Section~\ref{hyperbsection}). According to him, this idea is
originally due to Gabor Elek. I would like to thank Yaroslav
Kopylov for pointing out to me the reference \cite{GKS} where the
$L^p$-cohomology was first introduced. I am also grateful to Yves
de Cornulier, Pierre Pansu, Gilles Pisier, and Michael Puls for
their useful remarks and corrections.

\section{Main results}\label{mainresultsection}
{\it (The definitions of first $L^p$-cohomology, p-harmonic
functions and of first cohomology with values in a representation
are postponed to Section~\ref{prelimSection}.)}

Let $G$ be a locally compact group acting by measure-preserving
bijections on a measure space $(X,m)$. We say that the action is
strongly mixing (or mixing) if for every measurable subset of
finite measure $A\subset X$, $m(gA\cap A)\to 0$ when $g$ leaves
every compact subset of $G$. Let $\pi$ be the corresponding
continuous representation of $G$ in $L^p(X,m)$, where
$1<p<\infty$. In this paper, we will call such a representation a
mixing $L^p$-representation of $G$.

\begin{defn}\cite{T1}
Let $G$ be a locally compact, compactly generated group and let $S$
be a compact generating subset of $G$. We say that $G$ has Property
(CF) (Controlled F\o lner) if there exists a sequence of compact
subsets of positive measure $(F_n)$ satisfying the following
properties.
\begin{itemize}
\item $F_n\subset S^n$ for every $n$;

\item there is a constant $C<\infty$ such that for every $n$ and every $s\in S$,
$$\frac{\mu(sF_n\vartriangle F_n)}{\mu(F_n)}\leq C/n.$$
\end{itemize}
Such a sequence $F_n$ is called a controlled F\o lner\footnote{A
controlled F\o lner sequence is in particular a F\o lner sequence,
so that Property (CF) implies amenability.} sequence.
\end{defn}
In \cite{T1}, we proved that following family\footnote{This family
of groups also appears in \cite{coteva}.} of groups are (CF).
\begin{itemize}
\item[(1)] Polycyclic groups and connected amenable Lie groups;

\item[(2)] semidirect products
$\Z[\frac{1}{mn}]\rtimes_{\frac{m}{n}}\Z$, with $m,n$ co-prime
integers with $|mn|\ge 2$ (if $n=1$ this is the Baumslag-Solitar
group $BS(1,m)$); semidirect products $\left(\bigoplus_{i\in
I}\Q_{p_i}\right)\rtimes_{\frac{m}{n}}\Z$ with $m,n$ co-prime
integers, and $(p_i)_{i\in I}$  a finite family of primes (including
$\infty$: $\Q_{\infty}=\R)$) dividing $mn$;

\item[(3)] wreath products $F\wr\Z$ for $F$ a finite group.
\end{itemize}

Our main result is the following theorem.
\begin{ThmIntro}\label{Mainthm}
Let $G$ be a group with Property (CF) and let $\pi$ be a mixing
$L^p$-representation of $G$. Then the first reduced cohomology of
$G$ with values in $\pi$ vanishes, i.e. $\overline{H^1}(G,\pi)=0.$
\end{ThmIntro}

\noindent{\bf Invariance under quasi-isometry.} The proof of
\cite[Theorem~4.3.3]{Shal} that Property $H_{FD}$ is invariant under
quasi-isometry can be used identically in the context of
$L^p$-representations and replacing the hypothesis ``weak mixing" by
``mixing" since the induced representation of a mixing
$L^p$-representation is also a mixing $L^p$-representation. As a
result, we obtain that the property that $\overline{H^1}(G,\pi)=0$
for every mixing $L^p$-representation is invariant under
quasi-isometry between discrete amenable groups. It is also stable
by passing to (and inherited by) co-compact lattices in amenable
locally compact groups.

\medskip

It is well known \cite{Puls} that for finitely generated groups $G$,
the first reduced cohomology with values in the left regular
representation in $\ell^p(G)$ is isomorphic to the space $HD_p(G)$
of $p$-harmonic functions with gradient in $\ell^p$ modulo the
constants. We therefore obtain the following corollary.

\begin{CorIntro}\label{Liouville(L)}
Let $G$ be a discrete group with Property (CF). Then every
$p$-harmonic function on $G$ with gradient in $\ell^p$ is constant.
\end{CorIntro}
Using Von Neumann algebra technics, Cheeger and Gromov \cite{CG}
proved that every finitely generated amenable group $G$ has no
nonconstant harmonic function with gradient in $\ell^2$, the
generalization to every $1<p<\infty$ being conjectured by Gromov.

To obtain a version of Corollary \ref{Liouville(L)} for Lie
groups, we prove the following result (see
Theorem~\ref{Lp/actionprop}).

\begin{thm}\label{Lp/actionprop}
Let $G$ be a connected Lie group. Then for $1\leq p< \infty$, the
first $L^p$-cohomology of $G$ is topologically (canonically)
isomorphic to the first cohomology with values in the right
regular representation in $L^p(G)$, i.e.
$$H^1_p(G)\simeq H^1(G,\rho_{G,p}).$$
\end{thm}
Now, since this isomorphism induces a natural bijection
$$HD_p(G)\simeq \overline{H^1}(G,\rho_{G,p}),$$
we can state the following result that was conjectured by Pansu in
\cite{Pansu}. Recall that a Riemannian manifold is called closed at
infinity if there exists a sequence of compact subsets $A_n$ with
regular boundary $\partial A_n$ such that $\mu_{d-1}(\partial
A_n)/\mu_d(A_n)\to 0$, where $\mu_k$ denotes the Riemannian measure
on submanifolds of dimension $k$ of $M$.

\begin{CorIntro}\label{LiouvilleM}
Let $M$ be a homogeneous Riemannian manifold. If it is closed at
infinity, then for every $p>1$, every $p$-harmonic function on $M$
with gradient in $L^p(TM)$ is constant. In other words,
$HD_p(M)=0.$
\end{CorIntro}

Together with Pansu's results \cite[Théorème~1]{Pan'}, we obtain the
following dichotomy.

\begin{ThmIntro}\label{homogenThm}
Let $M$ be a homogeneous Riemannian manifold. Then the following
dichotomy holds.
\begin{itemize}
\item Either $M$ is quasi-isometric to a homogeneous Riemannian
manifold with strictly negative curvature, and then there exists
$p_0\geq 1$ such that $HD_p(M)\neq 0$ if and only if $p>p_0$;
\item or $HD_p(M)=0$ for every $p>1$.
\end{itemize}
\end{ThmIntro}

We also prove
\begin{ThmIntro}(see Corollary~\ref{hyperbthm})
A homogeneous Riemannian manifold $M$ has non-zero first reduced
$L^p$-cohomology for some $1<p<\infty$ if and only if it is
non-elementary\footnote{By non-elementary, we mean not
quasi-isometric to $\R$.} Gromov hyperbolic.
\end{ThmIntro}

To prove this corollary, we need to prove that a Gromov hyperbolic
Lie group has non-trivial first reduced $L^p$-cohomology for $p$
large enough. This is done in Section~\ref{hyperbsection}. Namely,
we prove a more general result.

\begin{ThmIntro}(see Theorem~\ref{hyperbolicLpthm})
Let $G$ be a Gromov hyperbolic metric measure space with bounded
geometry having a bi-Lipschitz embedded 3-regular tree, then for
$p$ large enough, it has non-trivial first reduced
$L^p$-cohomology.
\end{ThmIntro}

Corollary~\ref{hyperbthm} and Pansu's contribution to
Theorem~\ref{homogenThm} yield the following corollary.

\begin{CorIntro}
A non-elementary Gromov hyperbolic homogeneous Riemannian manifold
is quasi-isometric to a homogeneous Riemannian manifold with
strictly negative curvature.
\end{CorIntro}

(See \cite{Heintze} for an algebraic description of homogeneous
manifolds with strictly negative curvature).

\section{Organization of the paper.}

In the following section, we recall three definitions of first
cohomology:
\begin{itemize}
\item a coarse definition of the first $L^p$-cohomology on a
general metric measure space which is due to Pansu; \item the
usual definition of first $L^p$-cohomology on a Riemannian
manifold; \item the first cohomology with values in a
representation, which is defined for a locally compact group.
\end{itemize}
In Section~\ref{topoIsomsection}, we construct a natural
topological isomorphism between the $L^p$-cohomology of a
connected Lie group $G$ and the cohomology with values in the
right regular representation of $G$ in $L^p(G)$. We use this
isomorphism to deduce Corollary~\ref{LiouvilleM} from
Theorem~\ref{Mainthm}.

The proof of Theorem~\ref{Mainthm} splits into two steps. First
(see Theorem~\ref{sublinthm}), we prove that for any locally
compact compactly generated group $G$ and any mixing
$L^p$-representation $\pi$ of $G$, every 1-cocycle $b\in
Z^{1}(G,\pi)$ is {\it sublinear}, which means that for every
compact symmetric generating subset $S$ of $G$, we have
$$\|b(g)\|=o(|g|_S)$$
when $|g|_S\to \infty$, $|g|_S$ being the word length of $g$ with
respect to $S$. Then, we adapt to this context a remark that we
made with Cornulier and Valette (see~\cite[Proposition
3.6]{coteva}): for a group with Property (CF), a $1$-cocycle
belongs to $\overline{B}^1(G,\pi)$ if and only if it is sublinear.
The part ``only if" is an easy exercise and does not require
Property (CF). To prove the other implication, we consider the
affine action $\sigma$ of $G$ on $E$ associated to the $1$-cocycle
$b$ and use Property (CF) to construct a sequence of almost fixed
points for $\sigma$.

In Section~\ref{variante}, we propose a more direct
approach\footnote{However, the ingredients are the same:
sublinearity of cocycles, and existence of a controlled F\o lner
sequence.} to prove Corollary~\ref{LiouvilleM}. The interest is to
provide an explicit approximation of an element of $\D_p(G)$ by a
sequence of functions in $W^{1,p}(G)$ using a convolution-type
argument.

Finally, in Section~\ref{hyperbsection}, we prove that a Gromov
hyperbolic homogeneous manifold has non-trivial $L^p$-cohomology
for $p$ large enough. This section can be read independently.

\section{Preliminaries}\label{prelimSection}

\subsection{A coarse notion of first $L^p$-cohomology on a metric measure
space}\label{lpcoarseversion}

The following coarse notion of (first) $L^p$-cohomology is
essentially due to \cite{Pans} (see also the chapter about
$L^p$-cohomology in \cite{Gr}).

Let $X=(X,d,\mu)$ be a metric measure space, and let $p\geq 1$.
For all $s>0$, we write $\Delta_s=\{(x,y)\in X^2, d(x,y)\leq s\}.$

First, let us introduce the $p$-Dirichlet space $\D_p(X)$.
\begin{itemize}
\item The space $D_p(X)$ is the set of measurable functions $f$ on
$X$ such that
$$\int_{\Delta_s}|f(x)-f(y)|^pd\mu(x)d\mu(y)<\infty$$
for every $s>0$.

\item Let $\D_p(X)$ be the Banach space $D_p(X)/\C$ equipped with
the norm
$$\|f\|_{D_p}=\left(\int_{\Delta_1}|f(x)-f(y)|^p\mu(x)d\mu(y)\right)^{1/p}.$$
\item By a slight abuse of notation, we identify $L^p$ with its
image in $\D_p$.
\end{itemize}
\begin{defn}
The first $L^p$-cohomology of $X$ is the space
$$H^1_p(X)=\D_p(X)/L^p(X),$$ and the first reduced $L^p$-cohomology
of $X$ is the space
$$\overline{H^1}_p(X)=\D_p(X)/\overline{L^p(X)}^{\D_p(X)}.$$
\end{defn}

\begin{defn}{\bf ($1$-geodesic spaces)} We say that a metric space
$X=(X,d)$ is $1$-geodesic if for every two points $x,y\in X$,
there exists a sequence of points $x=x_1,\ldots x_m=y,$ satisfying
\begin{itemize}
\item $d(x,y)=d(x_1,x_2)+\ldots+d(x_{m-1},x_m)$,

\item for all $1\leq i\leq m-1$, $d(x_i,x_{i+1})\leq 1.$
\end{itemize}
\end{defn}

\begin{rem}
Let $X$ and $Y$ be two $1$-geodesic metric measure spaces with
bounded geometry in the sense of \cite{Pans}. Then it follows from
\cite{Pans} that if $X$ and $Y$ are quasi-isometric, then
$H^1_p(X)\simeq H^1_p(Y)$ and $\overline{H^1}_p(X)\simeq
\overline{H^1}_p(Y).$
\end{rem}

\begin{ex}
Let $G$ be a locally compact compactly generated group, and let
$S$ be a symmetric compact generating set. Then the word metric on
$G$ associated to $S$, $$d_S(g,h)=\in\{n\in\N, g^{-1}h\in S^n\},$$
defines a $1$-geodesic left-invariant metric on $G$. Moreover, one
checks easily two such metrics (associated to different $S$) are
bilipschitz equivalent. Hence, by Pansu's result, the first
$L^p$-cohomology of $(G,\mu, d_S)$ does not depend on the choice
of $S$.
\end{ex}

\begin{defn}{\bf (coarse notion of $p$-harmonic functions)}
Let $f\in D_p(X)$ and assume that $p>1$. The
$p$-Laplacian\footnote{Here we define a coarse $p$-Laplacian at
scale $1$: see \cite[Section~2.2]{T2} for a more general
definition.} of $f$ is
$$\Delta_p f(x)=\frac{1}{V(x,1)}\int_{d(x,y)\leq 1}|f(x)-f(y)|^{p-2}(f(x)-f(y))d\mu(y),$$
where $V(x,1)$ is the volume of the closed ball $B(x,1).$ A
function $f\in D_p(X)$ is called $p$-harmonic if $\Delta_p f=0$.
Equivalently, the $p$-harmonic functions are the minimizers of the
variational integral
$$\int_{\Delta_1}|f(x)-f(y)|^{p}d\mu(x)d\mu(y).$$
\end{defn}

\begin{defn}
We say that $X$ satisfies a Liouville $D_p$-Property if every
$p$-harmonic function on $X$ is constant.
\end{defn}

As $\D_p(X)$ is a strictly convex, reflexive Banach space, every
$f\in \D_p(X)$ admits a unique projection $\tilde{f}$ on the
closed subspace $\overline{L^p(X)}$ such that
$d(f,\tilde{f})=d(f,\overline{L^p(X)}).$ One can easily check that
$f-\tilde{f}$ is $p$-harmonic. In conclusion, the reduced
cohomology class of $f\in \D_p(X)$ admits a unique $p$-harmonic
representant modulo the constants. We therefore obtain
\begin{prop}
A metric measure space $X$ has Liouville $D_p$-Property if and
only if $\overline{H_p}^{1}(X)=0.$
\end{prop}

\subsection{First $L^p$-cohomology on a Riemannian manifold}\label{lpRiemsection}

Let $M$ be Riemannian manifold, equipped with its Riemannian
measure $m$. Let $1\leq p<\infty.$

Let us first define, in this differentiable context, the
$p$-Dirichlet space $\D_p$.

\begin{itemize}
\item Let $D_p$ be the vector space of continuous functions whose
gradient is (in the sense of distributions) in $L^p(TM)$.

\item Equip $D_p(M)$ with a pseudo-norm $\|f\|_{D_p}=\|\nabla
f\|_p$, which induces a norm on $D_p(M)$ modulo the constants.
Denote by $\D_p(M)$ the completion of this normed vector space.
\item Write $W^{1,p}(M)=L^p(M)\cap D_p(M)$. By a slight abuse of
notation, we identify $W^{1,p}(M)$ with its image in $\D_p(M)$.
\end{itemize}

\begin{defn}
The first $L^p$-cohomology of $M$ is the quotient space
$$H_p^1(M)=\D_p(M)/W^{1,p}(M),$$ and the first reduced $L^p$-cohomology
of $M$  is the quotient
$$\overline{H_p}^1(M)=\D_p(M)/\overline{W^{1,p}(M)},$$ where
$\overline{W^{1,p}(M)}$ is the closure of $W^{1,p}(M)$ in the
Banach space $\D_p(M)$.
\end{defn}
\begin{defn}{\bf ($p$-harmonic functions)}
A function $f\in D_p(M)$ is called $p$-harmonic if it is a weak
solution of
$$\textnormal{div}(|\nabla f|^{p-2}\nabla f)=0,$$
that is,
$$\int_M\langle|\nabla f|^{p-2}\nabla f,\nabla \varphi\rangle dm=0,$$
for every $\varphi\in C_0^{\infty}.$ Equivalently, $p$-harmonic
functions are the minimizers of the variational integral
$$\int_{M}|\nabla f|^{p}dm.$$
\end{defn}

\begin{defn}
We say that $M$ satisfies a Liouville $D_p$-Property if every
$p$-harmonic function on $M$ is constant.
\end{defn}

As $\D_p(M)$ is a strictly convex, reflexive Banach space, every
$f\in \D_p(M)$ admits a unique projection $\tilde{f}$ on the
closed subspace $\overline{W^{1,p}(M)}$ such that
$d(f,\tilde{f})=d(f,\overline{W^{1,p}(M)}).$ One can easily check
that $f-\tilde{f}$ is $p$-harmonic. In conclusion, the reduced
cohomology class of $f\in \D_p(M)$ admits a unique $p$-harmonic
representant modulo the constants. Hence, we get the following
well-known fact.

\begin{prop}
A Riemannian manifold $M$ has Liouville $D_p$-Property if and only
if $\overline{H_p}^{1}(M)=0.$
\end{prop}

\begin{rem}\label{pansuEqLprem}
In \cite{Pans}, Pansu proves (in particular) that if a Riemannian
manifold has bounded geometry (which is satisfied by a homogeneous
manifold), then the first $L^p$-cohomology defined as above is
topologically isomorphic to its coarse version defined at the
previous section. In particular, the Liouville $D_p$-Property is
invariant under quasi-isometry between Riemannian manifolds with
bounded geometry.
\end{rem}

\subsection{First cohomology with values in a representation}

Let $G$ be a locally compact group, and $\pi$ a continuous linear
representation  on a Banach space $E=E_{\pi}$. The space
$Z^1(G,\pi)$  is defined as the set of continuous functions
$b:G\to E$ satisfying, for all $g,h$ in $G$, the 1-cocycle
condition $b(gh)=\pi(g)b(h)+b(g)$. Observe that, given a
continuous function $b:G\to E$, the condition $b\in Z^1(G,\pi)$ is
equivalent to saying that $G$ acts by affine transformations on
$E$ by $\alpha(g)v=\pi(g)v+b(g)$. The space $Z^1(G,\pi)$ is
endowed with the topology of uniform convergence on compact
subsets.

The subspace of coboundaries $B^1(G,\pi)$ is the subspace (not
necessarily closed) of $Z^1(G,\pi)$ consisting of functions of the
form $g\mapsto v-\pi(g)v$ for some $v\in E$. In terms of affine
actions, $B^1(G,\pi)$ is the subspace of affine actions fixing a
point.

The first cohomology space of $\pi$ is defined as the quotient
space
$$H^1(G,\pi)=Z^1(G,\pi)/B^1(G,\pi).$$


The first {\it reduced} cohomology space of $\pi$ is defined as
the quotient space
$$\overline{H^1}(G,\pi)=Z^1(G,\pi)/\overline{B^1}(G,\pi),$$
where $\overline{B^1}(G,\pi)$ is the closure of $B^1(G,\pi)$ in
$Z^1(G,\pi)$ for the topology of uniform convergence on compact
subsets. In terms of affine actions, $\overline{B^1}(G,\pi)$ is
the space of actions $\sigma$ having almost fixed points, i.e. for
every $\eps>0$ and every compact subset $K$ of $G$, there exists a
vector $v\in E$ such that for every $g\in K$,
$$\|\sigma(g)v-v\|\leq \eps.$$
If $G$ is compactly generated and if $S$ is a compact generating
set, then this is equivalent to the existence of a sequence of
almost fixed points, i.e. a sequence $v_n$ of vectors satisfying
$$\lim_{n\to\infty}\sup_{s\in S}\|\sigma(s)v_n-v_n\|=0.$$

\section{$L^p$-cohomology and affine actions on $L^p(G)$.}\label{topoIsomsection}

Let $G$ be a locally compact group equipped with a left-invariant
Haar measure. Let $G$ act on $L^p(G)$ by right translations, which
defines a representation $\rho_{G,p}$ defined by
$$\rho_{G,p}(g)f(x)=f(xg)\quad \forall f\in L^p(G).$$
Note that this representation is isometric if and only if $G$ is
unimodular, in which case $\rho_{G,p}$ is isomorphic to the left
regular representation $\lambda_{G,p}$. In particular, in this
case, we can identify the first reduced cohomologies.

Now suppose that the group $G$ is also compactly generated and
equipped with a word metric $d_S$ associated to a compact
symmetric generating subset $S$. In this section, we prove that
the first cohomology with values in the regular
$L^p$-representation $\rho_{G,p}$ is topologically isomorphic to
the first $L^p$-cohomology $H_p^{1}(G)$ (here, we mean the coarse
version, see Section~\ref{lpcoarseversion}). By the result of
Pansu mentioned in Remark~\ref{pansuEqLprem}, if $G$ is a
connected Lie group equipped with left-invariant Riemannian metric
$m$, we can also identify $H^1(G,\rho_{G,p})$ with the first
$L^p$-cohomology on $(G,m)$ (see Section~\ref{lpRiemsection}). We
also obtain a direct proof of this fact.

\

We consider here the two following contexts: where $G$ is a
compactly generated locally compact group equipped with a length
function $d_S$; or $G$ is a connected Lie group, equipped with a
left-invariant Riemannian metric.

Consider the linear map $J:\; \D_p(G)\to Z^1(G,\rho_{G,p})$
defined by
$$J(f)(g)=b(g)=f-\rho_{G,p}(g)f.$$ $J$ is clearly well defined and
induces a linear map $HJ: \; H^1_p(G)\to H^1(G,\rho_{G,p}).$
\begin{thm}\label{Lp/actionprop}
For $1\leq p< \infty$, the canonical map $HJ:\;H^1_p(G)\to
H^1(G,\rho_{G,p})$ is an isomorphism of topological vector spaces.
\end{thm}
Let us start with a lemma.

\begin{lem}\label{lemcocycleregularise}
Let $1\leq p<\infty$ and $b\in Z^1(G,\rho_{G,p})$. Then there
exists a $1$-cocycle $c$ in the cohomology class of $b$ such that
\begin{enumerate}
\item the map $(g,x)\mapsto c(g)(x)$ is continuous;\item the
continuous map $f(x)=c(x^{-1})(x)$ satisfies
$c(g)=f-\rho_{G,p}(g)f$;

\item moreover if $G$ is a Riemannian connected Lie group, then
$c$ can be chosen such that $f$ lies in $D_p(G)$ (and in
$C^{\infty}(G)$).
\end{enumerate}
\end{lem}

\noindent{\bf Proof of the lemma.}  Let $\psi$ be a continuous,
compactly supported probability density on $G$. We define $c\in
Z^1(G,\rho_{G,p})$ by
$$c(g)=\int_G b(gh)\psi(h)dh-\int_Gb(h)\psi(h)dh=\int_G b(h)(\psi(g^{-1}h)-\psi(h))dh.$$
Clearly it satisfies $c(1)=0$ and we have
\begin{eqnarray*}
c(gg') & = & \int_G b(gg'h)\psi(h)dh-\int_G b(h)\psi(h)dh\\
       & = & \rho_{G,p}(g)\int_G b(g'h)\psi(h)dh+ \int_G b(g)\psi(h)dh-\int_G b(h)\psi(h)dh
\end{eqnarray*}
But note that
\begin{eqnarray*}
\int_G b(g)\psi(h)dh & = & \int_G b(ghh^{-1})\psi(h)dh\\
       & = & \rho_{G,p}(g)\int_G \rho_{G,p}(h)b(h^{-1})\psi(h)dh+ \int_G
       b(gh)\psi(h)dh\\
       & = & -\rho_{G,p}(g)\int_Gb(h)\psi(h)dh+ \int_G
       b(gh)\psi(h)dh.
\end{eqnarray*}
So we obtain
\begin{eqnarray*}
c(gg') & = & \rho_{G,p}(g)\left(\int_G b(g'h)\psi(h)dh-\int_G
b(h)\psi(h)dh\right)+\int_G
       b(gh)\psi(h)dh-\int_G b(h)\psi(h)dh\\
       & = & \rho_{G,p}(g)c(g')+c(g).
\end{eqnarray*}
So $c$ is a cocycle.

Let us check that $c$ belongs to the cohomology class of $b$.
Using the cocycle relation, we have
\begin{eqnarray*}
c(g)&=&\int_G (\rho_{G,p}(g)b(h)-b(g)\psi(h)dh-\int_Gb(h)\psi(h)dh\\
&=& b(g)+\int_G (\rho_{G,p}(g)b(h)-b(h)\psi(h)dh)\\
&=& b(g)+\rho_{G,p}(g)\int_G b(h)\psi(h)dh-\int_G b(h)\psi(h)dh.
\end{eqnarray*}
But since $\int_G b(h)\psi(h)dh\in L^p(G)$, we deduce that $c$
belongs to the cohomology class of $b.$

Now, let us prove that $(g,x)\mapsto c(g)(x)$ is continuous. It is
easy to see from the definition of $c$ that $g\mapsto c(g)(x)$ is
defined and continuous for almost every $x$: fix such a point
$x_0$. We conclude remarking that the cocycle relation implies
$$c(g)(x)=c(xg)(x_0)-c(g)(x_0).$$
Now we can define $f(x)=c(x^{-1})(x)$ and again the cocycle
relation for $b$ implies that $c(g)=f-\rho_{G,p}(g)f$.

Finally, assume that $G$ is a Lie group and choose a smooth
$\psi$. The function $\hat{\psi}$ defined by
$$\hat{\psi}(g)=\psi(g^{-1})$$ is also smooth and compactly
supported. We have
$$c(g)(x)=f(x)-f(xg)=\int_G b(h)(x)(\hat{\psi}(h^{-1}g)-\hat{\psi}(h^{-1}))dh.$$
Hence, $f$ is differentiable and
$$\nabla f (x)=\int_G b(h)(x)(\nabla \hat{\psi})(h^{-1}))dh,$$
and so $\nabla f\in L^p(TG).$ \epr

\medskip

\noindent{\bf Proof of Theorem~\ref{Lp/actionprop}.}  The last
statement of the lemma implies that $HJ$ is surjective. The
injectivity follows immediately from the fact that $f$ is
determined up to a constant by its associated cocycle $b=I(f)$.

We now have to prove that the isomorphism $HJ$ is a topological
isomorphism. This is immediate in the context of the coarse
$L^p$-cohomology. Let us prove it for a Riemannian connected Lie
groups. Let $S$ be a compact generating subset of $G$ and define a
norm on $Z^1(G,\rho_{G,p})$ by
$$\|b\|=\sup_{s\in S}\|b(s)\|_p.$$

Let $\psi$ be a regular, compactly generated probability density
on $G$ as in the proof of Lemma~\ref{lemcocycleregularise}. Denote
$$f\ast \psi(x)=\int_G f(xh)\psi(h).$$
We have
\begin{lem}\label{lemtopoisom}
There exists a constant $C<\infty$ such that for every $f\in
\D_p(G)$,
$$C^{-1}\|f\ast \psi\|_{\D_p}\leq \|J(f)\|\leq C\|f\|_{\D_p}.$$
\end{lem}
\noindent{\bf Proof of the lemma.}  First, one checks easily that
if $b$ is the cocycle associated to $f$, then the regularized
cocycle $c$ constructed in the proof of
Lemma~\ref{lemcocycleregularise} is associated to $f\ast \psi.$

We have
\begin{eqnarray*}
 \nabla (f\ast \psi)(x) & = & \int(f(xh)-f(x))\nabla_x(\psi(xh))dh\\
                         & = & \int(f(xh)-f(x))\nabla_x(\hat{\psi}(h^{-1}x))dh
\end{eqnarray*}
So
\begin{eqnarray*}
\|\nabla(f\ast \psi)\|_p & \leq & \sup_{h\in \supp(\hat{\psi})}\int|f(xh)-f(x)|^p\|\nabla\hat{\psi}\|_{\infty}^pdx\\
                         & = & \sup_{h\in
                         \supp(\hat{\psi})}\|b(h)\|^p\|\nabla\hat{\psi}\|_{\infty}^pdx,
\end{eqnarray*}
which proves the left-hand inequality of Lemma~\ref{lemtopoisom}.
Let $g\in G$ and $\gamma:\;[0,d(1,g)]\to G$ be a geodesic between
$1$ and $g$. For any $f\in \D_p(G)$ and $x\in G$, we have
$$(f-\rho_{G,p}(g)f)(x)=f(x)-f(xg)=\int_0^{d(1,g)}\nabla f(x)\cdot\gamma'(t)dt.$$
So we deduce that
$$\|f-\rho_{G,p}(g)f\|_p\leq d(1,g)\|\nabla f\|_p,$$
which proves the right-hand inequality of
Lemma~\ref{lemtopoisom}.\epr

\medskip

Continuity of $HJ$ follows from continuity of $J$ which is an
immediate consequence of Lemma~\ref{lemtopoisom}.

Let us prove that the inverse of $HJ$ is continuous. Let $b_n$ be
a sequence in $Z^1(G,\rho_{G,p})$, converging to $0$ modulo
$B^1(G,\rho_{G,p})$. This means that there exists a sequence $a_n$
in $B^1(G,\rho_{G,p})$ such that $\|b_n+a_n\|\to 0$. By
Lemma~\ref{lemcocycleregularise}, we can assume that
$b_n(g)=f_n-\rho_{G,p}(g)f_n$ with $f\in \D_p(G).$ On the other
hand, $a_n=h-\rho_{G,p}(g)h$ with $h\in L^p(G)$. As compactly
supported, regular\footnote{Regular here, means either continuous,
or $C^{\infty}$ if $G$ is a Lie group.} functions on $G$ are dense
in $L^p(G)$, we can assume that $h$ is regular. So finally,
replacing $f_n$ by $f_n+h_n$, which is in $\D_p(G)$, we can assume
that $J(f_n)\to 0$. Then, by Lemma~\ref{lemtopoisom},
$\|f_n\ast\psi\|_{\D_p}\to 0$. But by the proof of
Lemma~\ref{lemcocycleregularise}, $f_n\ast\psi$ is in the class of
$L^p$-cohomology of $f_n$. This finishes the proof of
Theorem~\ref{Lp/actionprop}. \epr

\section{Sublinearity of cocycles}

\begin{thm}\label{sublinthm}
Let $G$ be a locally compact compactly generated group and let $S$
be a compact symmetric generating subset. Let $\pi$ be a mixing
$L^p$-representation of $G$. Then, every $1$-cocycle $b\in
Z^1(G,\pi)$ is sublinear, i.e.
$$\|b(g)\|=o(|g|_S)$$
when $|g|_S\to \infty$, $|g|_S$ being the word length of $g$ with
respect to $S$.
\end{thm}
Let $L^p(X,m)$ the $L^p$-space on which $G$ acts. We will need the
following lemma.
\begin{lem}\label{lem1}
Let us keep the assumptions of the theorem. There exists a constant
$C<\infty$ such that for any fixed $j\in \N$,
$$\|\pi(g_1)v_1+\ldots +\pi(g_j)v_j\|_p^p\to \|v_1\|_p^p+\ldots +\|v_j\|_p^p$$
when $d_S(g_k,g_l)\to \infty$ whenever $k\neq l$, uniformly with
respect to $(v_1,\ldots,v_j)$ on every compact subset of
$(L^p(X,m))^j$.
\end{lem}

\noindent{\bf Proof of Lemma~\ref{lem1}.}
First, let us prove that if the lemma holds pointwise with respect to $\overline{v}=(v_1,\ldots,v_j)$, then it holds uniformly on every compact subset $K$ of $(L^p(X,m))^j$. Let us fix some $\eps>0$. Equip $(L^p(X,m))^j$ with the norm $$\|\overline{v}\|=\max_{i}\|v_i\|_p,$$ and take a finite covering of $K$ by balls of radius $\eps$: $B(\overline{w},\eps)$, $\overline{w}\in W$, where $W$ is a finite subset of $K$. Take $\min_{1\leq k\neq i\leq j}d_S(g_k,g_l)$ large enough so that for any $\overline{w}\in W$, $\|\pi(g_1)v_1+\ldots +\pi(g_j)v_j\|_p^p$ is closed to $\|v_1\|_p^p+\ldots +\|v_j\|_p^p$ up to $\eps$. As $\pi(g)$ preserves the $L^p$-norm for every $g\in G$, we immediately see that for any $\overline{v}$ in $K$, $\|\pi(g_1)v_1+\ldots +\pi(g_j)v_j\|_p^p$ is closed to $\|v_1\|_p^p+\ldots +\|v_j\|_p^p$ up to some $\eps'$ only depending on $K$, $p$ and $\eps$, and such that $\eps'\to 0$ when $\eps\to 0$.

So now, we just have to
prove the lemma for $v_1,\ldots,v_j$ belonging to a dense subset of $L^p(X,m)$.
Thus, assume that for every $1\leq k\leq j$, $v_k$ is bounded and compactly supported.
Let us denote by $A_k$ the support of $v_k$.
For every finite sequence $\overline{g}=g_1,\ldots, g_j$ of elements in $G$, we write, for every $1\leq i\leq j,$
\begin{itemize}
\item $U_{i,\overline{g}}=\left(\bigcup_{l\neq i}g_l A_l\right)\cap g_i
A_i;$

\item $A_{i,\overline{g}}=g_i  A_i\smallsetminus U_{i,\overline{g}}.$
\end{itemize}

The key point of the proof is  the following observation
\begin{clai}
For every $1\leq
i\leq j$,
$$m(U_{i,\overline{g}})\to 0,$$
when the relative distance between the $g_k$ goes to $\infty$.
\end{clai}
\noindent{\bf Proof of the claim.} For $u,v\in L^2(G,m)$,
write $\langle u,v\rangle =\int_{X}u(x)v(x) dm(x)$. For every $1\leq
i\leq j$,
\begin{eqnarray*}
m\left(\left(\bigcup_{l\neq i}g_l A_l\right)\cap g_i
A_i\right)& = & \left\langle
\sum_{l\neq i}\pi(g_l)1_{A_l},\pi(g_i)1_{A_i}\right\rangle\\
                                              & = & \sum_{l\neq
                                              i}\left\langle
                                              \pi(g_l)1_{A_l},\pi(g_i)1_{A_i}\right\rangle\\
                                              & = & \sum_{l\neq
                                              i}\left\langle
                                              \pi(g_l^{-1}g_i)1_{A_i},1_{A_l}\right\rangle\\
                                              & = & \sum_{l\neq
                                              i} m(g_l^{-1}g_i A_i\cap A_l)\to
                                              0
\end{eqnarray*}
by mixing property of the action.  \epr

\

\noindent{\bf Proof of the lemma.}
First, observe that by the claim,
    $$\|\pi(g_i)v_i1_{U_{i,\bar{g}}}\|_p^p\leq \|v_i\|_{\infty}^pm(U_{i,\overline{g}})\to 0,$$  when the relative distance between the $g_k$ goes to $\infty$. In other words, as $\pi(g_i)v_i=\pi(g_i)v_i1_{A_{i,\bar{g}}}+\pi(g_i)v_i1_{U_{i,\bar{g}}}$,
$$\|\pi(g_i)v_i1_{A_{i,\bar{g}}}-\pi(g_i)v_i\|_p^p\to 0.$$
In particular,
$$\|\pi(g_i)v_i1_{A_{i,\bar{g}}}\|_p^p\to\|v_i\|_p^p.$$
On the other hand, the $A_{i,\bar{g}}$ are piecewise disjoint. So finally, we have
\begin{eqnarray*}
\lim_{d_S(g_l,g_k)\to \infty}\|\pi(g_1)v_1+\ldots +\pi(g_j)v_j\|_p^p & =& \lim_{d_S(g_l,g_k)\to \infty}\|\pi(g_1)v_11_{A_{1,\overline{g}}}+\ldots +\pi(g_j)v_j1_{A_{j,\overline{g}}}\|_p^p\\
&=&  \lim_{d_S(g_l,g_k)\to \infty}\|\pi(g_1)v_11_{A_{1,\overline{g}}}\|^p+\ldots +\|\pi(g_j)v_j1_{A_{j,\overline{g}}}\|_p^p\\ &=&  \|v_1\|_p^p+\ldots +\|v_j\|_p^p,
\end{eqnarray*}
which proves the lemma. \epr

\

\noindent{\bf Proof of Theorem~\ref{sublinthm}.} Fix some $\eps>0$.
Let $g=s_1\ldots s_n$ be a minimal decomposition of $g$ into a
product of elements of $S$. Let $m\leq n$, $q$ and $r<m$ be positive
integers such that $n=qm+r$. To simplify notation, we assume $r=1$.
For $1\leq i< j\leq n$, denote by $g_j$ the prefix $s_1\ldots s_j$
of $g$ and by $g_{i,j}$ the subword $s_{i+1}\ldots s_j$ of $g$.
Developing $b(g)$ with respect to the cocycle relation, we obtain
$$b(g)=b(s_1)+\pi(g_1)b(s_2)+\ldots +\pi(g_{n-1})b(s_n).$$
Let us put together the terms in the following way
\begin{eqnarray*}
b(g)=  \left[b(s_1)+\pi(g_m)b(s_{m+1})+\ldots+ \pi(g_{(q-1)m})b(s_{(q-1)m+1})\right]\\
+\left[\pi(g_1)b(s_2)+\pi(g_{m+1})b(s_{m+2})+\ldots+ \pi(g_{(q-1)m+1})b(s_{(q-1)m+2})\right]\\
+\ldots+ \left[\pi(g_{m-1})b(s_m)+\pi(g_{2m-1})b(s_{2m})+\ldots+
\pi(g_{qm})b(s_{qm+1})\right]
\end{eqnarray*}
In the above decomposition of $b(g)$, consider each term between
$[\cdot]$, e.g. of the form
\begin{equation}\label{eq1}
\pi(g_k)b(s_{k+1})+\ldots +\pi(g_{(q-1)m+k})b(s_{(q-1)m+k+1})
\end{equation}
for $0\leq k\leq m-1$ (we decide that $s_0=1$). Note that since $S$
is compact and $\pi$ is continuous, there exists a compact subset
$K$ of $E$ containing $b(s)$ for every $s\in S$. Clearly since
$g=s_1\ldots s_n$ is a minimal decomposition of $g$, the length of
$g_{i,j}$ with respect to $S$ is equal to $j-i-1$. For $0\leq
i<j\leq q-1$ we  have
$$d_S(g_{im+k},g_{jm+k})=|g_{im+k,jm+k}|_S=(j-i)m\geq m.$$
So by Lemma~\ref{lem1}, for $m=m(q)$ large enough, the $p$-power of
the norm of (\ref{eq1}) is less than
$$\|b(s_{k+1})\|_p^p+\|b(s_{m+k+1})\|_p^p+\ldots +\|b(s_{(q-1)m+k+1})\|_p^p+1.$$
The above term is therefore less than $2q$. Hence, we have
$$\|b(g)\|_p\leq 2mq^{1/p}.$$
So for $q\geq q_0=(2/\eps)^{p/(p-1)}$, we have
$$\|b(g)\|_p/n\leq 2q^{1-1/p}\leq \eps.$$
Now, let $n$ be larger than $m(q_0)q_0$. We have $\|b(g)\|_p/|g|\leq
\eps.$\epr

\section{Proof of Theorem~\ref{Mainthm}}

Theorem~\ref{Mainthm} results from Theorem~\ref{sublinthm} and the
following result, which is an immediate generalization of
\cite[Proposition~3.6]{coteva}. For the convenience of the reader,
we give its short proof.

\begin{prop}\label{cotevathm}
Let $G$ be a group with property (CF) and let $\pi$ be a continuous
isometric action of $G$ on a Banach space $E$. Let $b$ a $1$-cocycle
in $Z^1(G,\pi)$. Then $b$ belongs to $\overline{B^1(G,\pi)}$ if and
only if $b$ is sublinear.
\end{prop}
\bpr Assume that $b$ is sublinear.

Let $(F_n)$ be a controlled F\o lner sequence in $G$. Define a
sequence $(v_n)\in E^{\N}$ by
$$v_n=\frac{1}{\mu(F_n)}\int_{F_n}b(g)dg.$$
We claim that $(v_n)$ defines a sequence of almost fixed points for
the affine action $\sigma$ defined by $\sigma(g)v=\pi(g)v+b(g).$
Indeed, we have
\begin{eqnarray*}
\|\sigma(s)v_n-v_n\|& = &\left\|\frac{1}{\mu(F_n)}\int_{F_n}\sigma(s)b(g)dg-\frac{1}{\mu(F_n)}\int_{F_n}b(g)dg\right\|\\
                & =& \left\|\frac{1}{\mu(F_n)}\int_{F_n}b(sg)dg-\frac{1}{\mu(F_n)}\int_{F_n}b(g)dg\right\|\\
                & = &\left\|\frac{1}{\mu(F_n)}\int_{s^{-1}F_n}b(g)dg-\frac{1}{\mu(F_n)}\int_{F_n}b(g)dg\right\|\\
                & \leq & \frac{1}{\mu(F_n)}\int_{s^{-1}F_n\vartriangle
                F_n}\|b(g)\|dg.
\end{eqnarray*}
Since $F_n\subset S^n$, we obtain that
$$\|\sigma(s)v_n-v_n\| \leq  \frac{C}{n}\sup_{|g|_S\leq
n+1}\|b(g)\|$$ which converges to $0$. This proves the non-trivial
implication of Proposition~\ref{cotevathm}.\epr

\section{Liouville $D_p$-Properties: a direct approach.}\label{variante}

In this section, we propose a direct proof of
Corollary~\ref{LiouvilleM}. Instead of using Theorem~\ref{Mainthm}
and Theorem~\ref{Lp/actionprop}, we reformulate the proof, only
using Theorem~\ref{sublinthm} and \cite[Theorem~11]{T1}. The
interest is to provide an explicit approximation of an element of
$\D_p(G)$ by a sequence of functions in $W^{1,p}(G)$ using a
convolution-type argument. Since Liouville $D_p$-Property is
equivalent to the vanishing of $\overline{H_p}^1(G)$, we have to
show that for every $p$-Dirichlet function on $G$, there exists a
sequence of functions $(f_n)$ in $W^{1,p}(G)$ such that the
sequence $(\|\nabla(f-f_n)\|_p)$ converges to zero. Let $(F_n)$ be
a {\it right} controlled F\o lner sequence. By a standard
regularization argument, we can construct for every $n$, a smooth
$1$-Lipschitz function $\varphi_n$ such that
\begin{itemize}
\item $0\leq \varphi_n\leq 1$;

\item for every $x\in F_n,$ $\varphi_n(x)=1$;

\item for every $y$ at distance larger than $2$ from $F_n$,
$\varphi_n(y)=0.$
\end{itemize}
Denote by $F_n'=\{x\in G: \; d(x,F_n)\leq 2\}$. As $F_n$ is a
controlled F\o lner sequence, there exists a constant $C<\infty$
such that
$$\mu(F_n'\smallsetminus F_n)\leq C\mu(F'_n)/n$$ and
$$F'_n\subset B(1,Cn).$$
Define $$p_n=\frac{\varphi_n}{\int_G\varphi_nd\mu}.$$ Note that
$p_n$ is a probability density satisfying for every $x\in X$,
$$|\nabla p_n(x)|\leq \frac{1}{\mu(F_n)}.$$
For every $f\in D_p(G)$, write $P_nf(x)=\int_X
f(y)p_n(y^{-1}x)d\mu(y).$ As $G$ is unimodular,
$$P_nf(x)=\int_X f(yx^{-1})p_n(y^{-1})d\mu(y).$$
We claim that $P_nf-f$ is in $W^{1,p}.$ For every $g\in G$ and every
$f\in D_p$, we have
$$\|f-\rho(g)f\|_p\leq d(1,g)\|\nabla f\|_p.$$
Recall that the support of $p_n$ is included in $F'_n$ which itself
is included in $B(1,Cn).$ Thus, integrating the above inequality, we
get
$$\|f-P_nf\|_p\leq Cn\|\nabla f\|_p,$$
so $f-P_nf\in L^p(G).$

It remains to show that the sequence $(\|\nabla P_nf\|_p)$
converges to zero. We have
$$\nabla P_nf(x) = \int_{G}f(y)\nabla p_n(y^{-1}x)d\mu(y)$$
Since $\int_G \nabla pd\mu=0$, we get
\begin{eqnarray*}
\nabla P_nf(x) & = & \int_{G}(f(y)-f(x^{-1}))\nabla
p_n(y^{-1}x)d\mu(y)\\
 & = & \int_{G}(f(yx^{-1})-f(x^{-1}))\nabla
p_n(y^{-1})d\mu(y).
\end{eqnarray*}
Hence,
\begin{eqnarray*}
\|\nabla P_nf\|_p & \leq & \int_G \|\lambda(y)f-f\|_p |\nabla
p_n(y^{-1})|d\mu(y)\\
& \leq & \frac{1}{\mu(F_n)}\int_{F'_n\smallsetminus F_n}
\|\lambda(y)f-f\|_pd\mu(y)\\
&\leq & \frac{\mu(F'_n\smallsetminus F_n)}{F_n}\sup_{|g|\leq
Cn}\|b(g)\|_p\\
&\leq &\frac{C}{n}\sup_{|g|\leq Cn}\|b(g)\|_p
\end{eqnarray*}
where $b(g)=\lambda(g)f-f$. Note that $b\in
Z^{1}(G,\lambda_{G,p})$. Thus, by Theorem~\ref{sublinthm},
$$\|\nabla P_nf\|_p \to 0.$$
This completes the proof of Corollary~\ref{LiouvilleM}.\epr

\section{Non-vanishing of the first reduced $L^p$-cohomology on a non-elementary Gromov hyperbolic space.}\label{hyperbsection}

Let us start with a remark about first $L^p$-cohomology on a
metric measure space.

\begin{rem} {\bf (Coupling between $1$-cycles and $1$-cocycles)}
A $1$-chains on $(X,d,\mu)$ is a functions supported on
$\Delta_r=\{(x,y)\in X^2, d(x,y)\leq r\}$ for some $r>0.$ The
$L^p$-norm of a (measurable) $1$-chain $s$ is the norm
$$\left(\int_{X^2}|s(x,y)|^pd\mu(x)d\mu(y)\right)^{1/p}.$$
A $1$-chain $s$ is called a $1$-cycle if $s(x,y)=s(y,x).$

Given $f\in \D_p$, we define a $1$-cocycle associated to $f$ by
$c(x,y)=f(x)-f(y)$, for every $(x,y)\in X^2$. Let $s$ be a
$1$-cycle in $L^q$, with $1/p+1/q=1$. We can form a coupling
between $c$ and $s$
$$\langle c,s\rangle=\int_{X^2} c(x,y)s(x,y)d\mu(x)d\mu(y)=\int_{X^2} (f(x)-f(y))s(x,y)d\mu(x)d\mu(y).$$
Clearly, if $f\in L^p$, then as $s$ is a cycle, we have $\langle
c,s\rangle=0$. This is again true for $f$ in the closure of
$L^p(X)$ for the norm of $\D_p(X)$. Hence, to prove that a
$1$-cocycle $c$ is non-trivial in $\overline{H^1}_p(X)$, it is
enough to find a $1$-cycle in $L^q$ whose coupling with $c$ is
non-zero.
\end{rem}

The main result of this section is the following theorem.

\begin{thm}\label{hyperbolicLpthm}
Let $X$ be a Gromov hyperbolic $1$-geodesic metric measure space
with bounded geometry having a bi-Lipschitz embedded 3-regular
tree, then for $p$ large enough, it has non-trivial first reduced
$L^p$-cohomology.
\end{thm}
From this theorem, we will deduce

\begin{CorIntro}\label{hyperbthm}
A homogeneous Riemannian manifold $M$ has non-zero first reduced
$L^p$-cohomology for some $1<p<\infty$ if and only if it is
non-elementary Gromov hyperbolic.
\end{CorIntro}
\noindent{\bf Proof of Corollary~\ref{hyperbthm}.} By
Theorem~\ref{homogenThm}, if $M$ has non-zero
$\overline{H}_p^1(M)$ for some $1<p<\infty$, then being
quasi-isometric to a negatively curved homogeneous manifold, it is
non-elementary Gromov hyperbolic.

Conversely, let $M$ be a Gromov hyperbolic homogeneous manifold.
As $M$ is quasi-isometric to its isometry group $G$, which is a
Lie group with finitely many components, we can replace $M$ by
$G$, and assume that $G$ is connected. If $G$ has exponential
growth, then \cite[Corollary~1.3]{CT} it has a bi-Lipschitz
embedded 3-regular tree $T$, and hence
Theorem~\ref{hyperbolicLpthm} applies. Otherwise $G$ has
polynomial growth, and we conclude thanks to the following
classical fact.
\begin{prop}
A non-elementary Gromov-hyperbolic connected Lie group has exponential growth.
\end{prop}

\bpr Let $G$ be a connected Lie group with polynomial growth. By
\cite{Gui}, $G$ is quasi-isometric to a simply connected nilpotent
group $G$, whose asymptotic cone \cite{PanCarnot} is homeomorphic
to another (graded) simply connected nilpotent Lie group with same
dimension. Hence, unless $G$ is quasi-isometric to $\R$, the
asymptotic cone of $G$ has dimension larger or equal than $2$. But
\cite[page~37]{Gr} the asymptotic cone of a Gromov hyperbolic
space is an $\R$-tree, and therefore has topological dimension 1.
\epr

\

\noindent{\bf Proof of Theorem~\ref{hyperbolicLpthm}.} The proof
contains ideas that we found in \cite[page~258]{Gr}. Roughly
speaking, we start by considering a non-trivial cycle defined on a
bi-Lipschitz embedded 3-regular subtree $T$ of $X$. To construct a
$1$-cocycle which has non-trivial reduced cohomology, we take a
Lipschitz function $F$ defined on the boundary of $X$, such that
$F$ is non-constant in restriction to the boundary of the subtree
$T$. We then extend $F$ to a function defined $f$ on $X$ which
defines a $1$-cocyle in $D_p(X)$. Coupling this cocycle with our
cycle on $T$ proves its non-triviality in $\overline{H}_p^1(X)$.

\

\noindent{\bf Boundary at infinity of a hyperbolic space.} To
denote the distance between to points in $X$ or in its boundary,
we will use indifferently the notation $d(x,y)$, or the notation
of Gromov $|x-y|.$ Let us fix a point $o\in X$. We will denote
$|x|=|x-o|=d(x,o).$

Consider the Gromov boundary (see \cite[Chapter~1.8]{Grohyperb} or
\cite{GH}) of $X$, i.e. the set of geodesic rays issued from $o$
up to Hausdorff equivalence.

For $\eps$ small enough, there exists \cite{GH} a distance
$|\cdot|_{\eps}$ on $\partial_{\infty}X$, and $C<\infty$ such that
$$|u-v|_{\eps}\leq \limsup_{t\to \infty} e^{-\eps(v(t)|w(t))}\leq C|u-v|_{\eps}.$$
for all $v,w\in \partial_{\infty}X$, where $(\cdot|\cdot)$ denotes
the Gromov product, i.e.
$$(x|y)=\frac{1}{2}(|x|+|y|-|x-y|).$$

\

\noindent{\bf Reduction to graphs.} A $1$-geodesic metric measure
space with bounded geometry is trivially quasi-isometric to a
connected graph with bounded degree (take a maximal $1$-separated
net, and join its points which are at distance $1$ by an edge).
Hence, we can assume that $X$ is the set of vertices of a graph
with bounded degree.

\

\noindent{\bf A Lipschitz function on the boundary.} By
\cite[page~221]{Gr}, $T$ has a cycle which has a non-zero pairing
with every non-zero 1-cochain $c$ on $T$ supported on a single
edge. Hence, to prove that $\overline{H}_p^1(X)\neq 0$, it is
enough to find an element $\overline{c}$ in $D_p(X)$ whose
restriction to $T$ is zero everywhere but on $e$.

The inclusion of $T$ into $X$ being bi-Lipschitz, it induces a
homeomorphic inclusion of the boundary of $T$, which is a Cantor
set, into the boundary of $X$. We therefore identify
$\partial_{\infty}T$ with its image in $\partial_{\infty}X$. Now,
consider a partition of $\partial_{\infty}T$ into two clopen
non-empty subsets $O_1$ and $O_2$. As $O_1$ and $O_2$ are disjoint
compact subsets of $\partial_{\infty}X$, they are at positive
distance from one another. Hence, for $\delta>0$ small enough, the
$\delta$-neighborhoods $V_1$ and $V_2$ of respectively $O_1$ and
$O_2$ in $\partial_{\infty}X$ are disjoint. Denote by $e$ the edge
of $T$ joining the two complementary subtrees $T_1$ and $T_2$
whose boundaries are respectively $O_1$ and $O_2$.

Now, take a Lipschitz function $F$ on  $\partial_{\infty}X$ which
equals $0$ on $V_1$ and $1$ on $V_2$.

\

\noindent{\bf Extension of $F$ to all of $X$.}

Let us first assume that every point in $X$ is at bounded distance
from a geodesic ray issued from $o$.

Let us define a function $f$ on $X$: for every $x$ in
$X\smallsetminus \{o\}$, we denote element by $u_x$ a geodesic ray
issued from $o$ and passing at distance at most $C$ from $x$ ($C$
is a constant). Define
$$f(x)=F(u_x)\quad \forall x\in X\smallsetminus \{o\}.$$
Let us prove that for $p$ large enough, $f\in D_p(X).$ Take two
elements $x$ and $y$ in $X$ such that $|x-y|\leq 1$, we have
\begin{eqnarray*}
|u_x(t)-u_y(t)| & \leq & |u_x(t)-x|+|x-u_y(t)|\\
& \leq & |u_x(t)-x|+|y-u_y(t)|+|x-y|\\
& \leq & |u_x(t)-x|+|y-u_y(t)|+1.
\end{eqnarray*}
So for large $t$,
\begin{eqnarray*}
(u_x(t)|u_y(t)) & = &|u_x(t)|+|u_y(t)|-|u_x(t)-u_y(t)|\\
 & \geq & |u_x(t)|+|u_y(t)|-|u_x(t)-x|-|y-u_y(t)|-1\\
 & = & |x|+|y|-2C-1\\
 & \geq & 2|x|-2(C+1)
\end{eqnarray*}
Hence, since $F$ is $K$-Lipschitz,
\begin{eqnarray*}
|f(x)-f(y)|^p & = & |F(u_x)-F(u_y)|^p\\
            & \leq & K^p |u_x-u_y|^p\\
            & \leq  & K^p\limsup_{t\to \infty} e^{-p\eps(u_x(t)|u_y(t))}\\
            & \leq & K^pe^{-2p\eps|x|+2(C+1)p}
\end{eqnarray*}
On the other hand, as $\mu(B(o,|x|))\leq Ce^{\lambda |x|}$ for
some $\lambda$, if $2p\eps>\lambda$, then $f$ is in $D_p(X)$.

Now, let us consider the values of $f$ along $T$. To fix the
ideas, let us assume that $o$ is a vertex of $T$. For $i=1,2,$
take $x_i$ a vertex of $T_i$. Let $e_{x_i}$ be the edge whose one
extremity is $x_i$ and that separates $o$ and $x_i$. Let $T_{x_i}$
be the connected component of $T\smallsetminus\{e_{x_i}\}$
contained in $T_i$. Let $y$ be a vertex of $T_{x_i}$ and let $z$
be an element of $\partial_{\infty}T_{x_i}$. One the one hand $z$
can be interpreted as a geodesic $t$ issued from $o$ and passing
through $y$ in $T$, and on the other hand as a geodesic ray $v$
issued from $o$ in $X$. As $T$ is bi-Lipschitz embedded in $X$,
$v$ is a quasi-geodesic ray in $X$. Hence it stays at bounded
distance --say less than $C$-- from $v$. In particular, $v$ passes
at distance less than $C$ from $y$. So for $d(o,y)$ large enough,
$|u_y-v|\leq \delta$. Hence, choosing ${x_i}$ far enough from $o$
in $T_i$, we have that all $u_y$ where $y\in T_{x_i}$ belong to
$V_i$. Now, up to modify $T$, we can suppose that $x_1$ and $x_2$
are the end points of the edge $e$. Then, along $T$, $c$ takes
only two values: $0$ on $T_1=T_{x_1}$ and $1$ on $T_2=T_{x_2}.$
Hence its coupling with the cycle of \cite[page~221]{Gr} is
non-zero, which implies that $\overline{H}_p^1(X)\neq 0$.

\

\noindent{\bf Reduction to the case when every point in $X$ is at
bounded distance from a geodesic ray issued from $o$.}

Consider the graph $\tilde{X}$ obtained by gluing a copy of $\N$
to every vertex of $X$. Clearly $\tilde{X}$ satisfies that every
point in $X$ is at bounded distance from a geodesic ray issued
from $o$. Applying the above to $\tilde{X}$, we construct an
element $\tilde{f}$ in $D_p(\tilde{X})$ that has a non-trivial
coupling with the cycle that we considered on $T$. Clearly the
restriction $f$ of $\tilde{f}$ to $X$ also satisfies these
properties and hence defines a non-trivial cocycle in
$\overline{H}_p^1(X)$. \epr

\bigskip
\footnotesize

\noindent \noindent Romain Tessera\\
Department of mathematics, Vanderbilt University,\\ Stevenson
Center, Nashville, TN 37240 United,\\ E-mail:
\url{tessera@clipper.ens.fr}

\end{document}